\documentclass[12pt]{article}
\usepackage{amsmath,amsthm,verbatim,amssymb,amsfonts,amscd, graphics, mathrsfs,hyperref}
\theoremstyle{plain}
\newtheorem{theorem}{Theorem}
\newtheorem*{intro-theorem}{Theorem}

\newtheorem*{intro-corollary}{Corollary}

\theoremstyle{definition}

\theoremstyle{remark}

\newcommand{\Z}{\mathbb{Z}}

\renewcommand{\char}{\text{char}}

\newcommand{\mc}[1]{\mathcal{#1}}

\makeatletter
\def\blfootnote{\xdef\@thefnmark{}\@footnotetext}
\makeatother
\begin{document}

\title{Corrigendum to ``Intersection homology with field coefficients: $K$-Witt spaces and $K$-Witt bordism''}
\author{Greg Friedman\footnote{This work was partially supported by a grant from the Simons Foundation (\#209127 to Greg Friedman)} }

\date{August 14, 2012}

\maketitle

The author's paper \cite{GBF21} concerns $K$-Witt spaces and, in particular, a computation of the bordism theory of such spaces. However, there is an error in the computation of the coefficient groups in dimensions $4k+2$ when $\char(K)=2$. In this corrigendum, we state, as far as possible, the correct results. Details can be found in \cite{GBF34}. 

If we consider $K$-Witt spaces and $K$-Witt bordism using $K$-orientations, then for  $\char(K)=2$, this is unoriented bordism, which we denote $\mc N_*^{K-\text{Witt}}$. 

\begin{theorem}
For a field $K$ with $\char(K)=2$ and for\footnote{Since these are geometric bordism groups, they vanish in negative degree.} $i\geq 0$,
\begin{equation*}
\mc N_i^{K-\text{Witt}}\cong
\begin{cases}
\Z_2, &i\equiv 0\mod 2,\\
0,  &i\equiv 1\mod 2.
\end{cases}
\end{equation*}
\end{theorem}

This result is also provided without detailed  proof by Goresky in \cite[page 498]{Go84}. 

If we consider $K$-Witt spaces and $K$-Witt bordism using $\Z$-orientations, then we denote the bordism theory by $\Omega_{*}^{K-\text{Witt}}$. In this case, there remains one ambiguity in the computation, but we can show the following:

\begin{theorem}\label{T}
For a field $K$ with $\char(K)=2$ and $k\geq 0$,
\begin{enumerate}
\item $\Omega_{0}^{K-\text{Witt}}\cong\Z$,
\item $\Omega_{4k}^{K-\text{Witt}}\cong \Z_2$,
\item $\Omega_{4k+1}^{K-\text{Witt}}=\Omega_{4k+3}^{K-\text{Witt}}=0$,
\item \label{I} 
  Either
\begin{enumerate}
\item $\Omega_{4k+2}^{K-\text{Witt}}=0$ for all $k$, or
\item there exists some $N>0$ such that $\Omega_{4k+2}^{K-\text{Witt}}=0$ for all $k<N$ and $\Omega_{4k+2}^{K-\text{Witt}}\cong \Z_2$ for all $k\geq N$. 
\end{enumerate}
\end{enumerate}
\end{theorem}

See \cite{GBF34} for a discussion of the difficulty in  determining which case of item \eqref{I} of Theorem \ref{T} holds. 

Independent of the existence or value of $N$ in condition \eqref{I} of the theorem, the computations from \cite[Section 4.5]{GBF21} of $\Omega^{K-\text{Witt}}_*(\,\cdot\,)$ as a generalized homology theory on CW complexes  continue to hold and to imply that for $\char(K)=2$, $$\Omega^{K-\text{Witt}}_n(X)\cong \bigoplus_{r+s=n}H_r(X;\Omega^{K-\text{Witt}}_s).$$ Similarly,
$$\mc N^{K-\text{Witt}}_n(X)\cong \bigoplus_{r+s=n}H_r(X;\mc N^{K-\text{Witt}}_s).$$

\paragraph{Other minor errata.}
In \cite{GBF21} it should not be part of the definition of a $K$-Witt space that the space be irreducible as a pseudomanifold. However, as every $K$-Witt space of dimension $>0$ is bordant to an irreducible $K$-Witt space \cite[page 1099]{Si83}, this error does not affect the bordism group computations of \cite{GBF21}. Not every $0$-dimensional $K$-Witt space is bordant to an irreducible one, but the computation of  $\Omega_{0}^{K-\text{Witt}}$ reduces to the manifold theory and gives the result of  \cite{GBF21} if one removes irreducibility from the definition.

The argument  that $\Omega_{4k+2}^{K-\text{Witt}}=0$ given in \cite{GBF21} does not hold when $k=0$. However, in this dimension it is not difficult to prove the result directly; details are provided in \cite{GBF34}.  

\providecommand{\bysame}{\leavevmode\hbox to3em{\hrulefill}\thinspace}
\providecommand{\MR}{\relax\ifhmode\unskip\space\fi MR }
\providecommand{\MRhref}[2]{%
  \href{http://www.ams.org/mathscinet-getitem?mr=#1}{#2}
}
\providecommand{\href}[2]{#2}

\end{document}